\documentclass[11pt]{article}
\usepackage{amsmath,amsfonts,enumerate,amssymb,verbatim,latexsym,mathrsfs,cases,mathdots,colordvi,url,extarrows,graphicx}
\usepackage{cite} % use short line for continuous multiple references
\usepackage{relsize}
\usepackage{mathrsfs}
\usepackage{xcolor}
\usepackage{indentfirst}%shou hang suo jin
\usepackage{float}

\renewcommand{\Box}{\rule{2.2mm}{2.2mm}}

\newcommand{\dist}{{\rm dist}}

\def\beginproof{\par\noindent {\bf Proof.}\ \ }
\def\endproof{\hskip .5cm $\Box$ \vskip .5cm}

\begin{document}
\newtheorem{theorem}{Theorem}[section]
\newtheorem{definition}[theorem]{Definition}
\newtheorem{lemma}[theorem]{Lemma}
\newtheorem{assump}[theorem]{Assumption}
\newtheorem{proposition}[theorem]{Proposition}
\newtheorem{corollary}[theorem]{Corollary}
\newtheorem{example}[theorem]{Example}
\newtheorem{remark}[theorem]{Remark}
\newtheorem{refer}{Reference}
\begin{sloppypar}
\title{On second-order weak sharp minima of general nonconvex set-constrained optimization problems
\thanks{The alphabetical order of the paper indicates the equal contribution to the paper.}}
\author{Xiaoxiao Ma\thanks{Department of Mathematics and Statistics, University of Victoria, Canada.  Email: xiaoxiaoma9638@gmail.com. } \and Wei Ouyang\thanks{School of Mathematics, Yunnan Normal University, Kunming 650500, People's Republic of China; Yunnan Key Laboratory of Modern Analytical Mathematics and Applications, Kunming 650500, China. The research of this author was partially supported by the National Natural Science Foundation of the People's Republic of China [Grant 12261109] and the Basic Research Program of Yunnan Province
[Grant 202301AT070080]. Email: weiouyangxe@hotmail.com. } \and Jane J. Ye \thanks{Corresponding author. Department of Mathematics and Statistics, University of Victoria, Canada. The research of this author was partially
supported by NSERC. Email: janeye@uvic.ca.} \and Binbin Zhang\thanks{School of Science, Kunming University of Science and Technology, Kunming 650500, People's Republic of China.  Email: bbzhang@kust.edu.cn.}}
  \date{}
\maketitle
\begin{abstract}
%	In this paper, we  characterize  the local  second-order weak sharp minima for a very general classes of nonconvex  optimization problems. By employing the classical/lower generalized support functions, we develop some new  second-order optimality conditions. These conditions are derived using the corresponding asymptotic second-order tangent cone and the outer second-order tangent set.  Our  results extend the existing ones by removing  the need for convexity on the constraint set and/or the outer second-order tangent set, as well as the requirement of nonemptiness of the outer second-order tangent set.  Moreover in our sufficient condition, unlike the classical results we do not require restrictive assumptions such as the uniform second-order regularity of the constraint set and the property of uniform approximation of the critical cones.
	
{This paper explores local second-order weak sharp minima for a broad class of nonconvex optimization problems. We propose novel second-order optimality conditions formulated through the use of classical and lower generalized support functions. These results are based on   asymptotic second-order tangent cones and outer second-order tangent sets. Specifically, our findings eliminate the necessity of assuming convexity in the constraint set and/or the outer second-order tangent set, or the nonemptiness of the outer second-order tangent set. Furthermore, unlike traditional approaches, our sufficient conditions do not rely on strong assumptions such as the uniform second-order regularity of the constraint set and the property of uniform approximation of the critical cones.}

\vskip 10 true pt

%\noindent {\bf {\color{cyan}Communicated by Fabian Flores-Bazan.}}

\vskip 10 true pt
\noindent {\bf Key words.}\quad second-order weak sharp minima,  optimality condition, support function, the lower generalized support function, outer second-order tangent set, asymptotic second-order tangent cone
\vskip 10 true pt


\end{abstract}
\newpage
\section{Introduction}

This paper is devoted to the study of the local second-order weak sharp minima for a general set-constrained optimization problem of the following form
\begin{equation}\label{op}
\min f(x)\,\,\, {\rm s.t.}\,\,\, g(x)\in K,
\end{equation}
where $f:\mathbb{R}^n\rightarrow \mathbb{R}$ and $g:\mathbb{R}^n\rightarrow \mathbb{R}^m$ are twice continuously differentiable, and  $K$ is a closed subset of $\mathbb{R}^m$.

%The problem \eqref{op} is a very general non-convex optimization problem since even when the set $K$ is convex, the feasible region of the problem \eqref{op} is still non-convex
%unless the function $g$ is affine. When $K$ is a convex cone, the problem \eqref{op} is referred to  a  cone-constrained optimization problem (e.g., \cite{lhp12}) which includes the nonlinear program, the second-order cone program and the semidefinite cone program as special cases. The possibility of $K$ being nonconvex in the problem \eqref{op} extends a  cone-constrained optimization problem to include much more general problems and provides a generic framework for a number of important optimization problem classes such as  the mathematical program with equilibrium constraints, the mathematical program with second-order cone complementarity constraints and the mathematical program with semidefinite cone complementarity constraints (\cite{chent05,chaosye14,gfer14,gferyr20,yezhou16}).

The optimization problem \eqref{op} represents a broad class of nonconvex problems. Even when the set 
$K$ exhibits convexity, the feasible set can still be nonconvex. In cases where $K$ forms a convex cone, \eqref{op} is commonly termed a cone-constrained optimization problem (e.g., \cite{lhp12}). This category encompasses well-known problem types such as nonlinear programs, second-order cone programs, and semidefinite programs. Allowing $K$ to be nonconvex further generalizes \eqref{op}, enabling it to model a wider array of challenging optimization problems. Notable examples include mathematical programs with equilibrium constraints (MPECs) \cite{gfer14,gferyr20}, second-order cone complementarity constraints \cite{%chent05,
yezhou16}, and semidefinite cone complementarity constraints \cite{chaosye14}.

%$C(\bar x):=\{d\in T_{g^{-1}(K)}(\bar x)\mid\nabla f(\bar x)d\leq 0\}$ denote the critical cone at a feasible solution $\bar x$ unless otherwise specified.
%Let $S$ be a closed subset of the feasible set $\Phi$ such that $f$ is finite and constant on $S$, that is, $$S:=\{x\in \Phi \mid f(x)=f_0\}\quad \text{for some}\ f_0\in \mathbb{R}.$$

{Throughout the paper, we denote by  $\Phi:=\{x\in \mathbb{R}^n \mid g(x)\in K\}= g^{-1}(K)$ the feasible set of problem (\ref{op}) and  $S$  the set of optimal solutions of problem \eqref{op}.} A feasible solution $\bar x\in \Phi$ is said to be a local second-order weak sharp minimizer or satisfy the local quadratic growth condition with respect to the solution set $S$, cf. \cite[Definition 3.141]{shap} of problem \eqref{op}, provided  that there exist positive numbers $\kappa,\delta>0$ such that
\begin{equation}\label{los}
f(x)\geq f(\bar x) + \kappa [\dist (x,S)]^2,\,\,\,\forall x\in \Phi\cap B_{\delta}(\bar x),
\end{equation}
where  $B_{\delta}(\bar x)$ is the open ball centered at $\bar x$ with radius $\delta$.  If $\bar x$ is an isolated point of the solution set $S$, then condition \eqref{los} reduces to
the local quadratic growth condition at $\bar x$:
\begin{equation}\label{quadraticgrow}
f(x)\geq f(\bar x) + \kappa \|x-\bar x\|^2,\,\,\,\forall x\in \Phi\cap B_{\delta}(\bar x),\end{equation}
for a sufficiently small $\delta>0$, while in general the condition \eqref{los} is weaker than \eqref{quadraticgrow}.  Building on the work in \cite{OuyangYe}, where first- and second-order optimality conditions for a local minimizer of problem \eqref{op} were established under the assumption that the set $K$ may be nonconvex and without requiring the second-order tangent set to be convex and nonempty, we extend their framework to investigate local second-order weak sharp minima, a special subset of the set of local minimizers. Specifically, we consider this concept in the setting where $K$ is nonconvex and the second-order tangent set is not necessarily convex and nonempty.
%If the set $S$ is compact, then the quadratic growth condition holds at every point of $S$ if and only if the local quadratic growth condition with respect to the solution set $S$.
%Since  the solution set of a convex optimization problem must be convex, it is likely that there is no isolated minima for a convex problem unless the set of solution is non-unique.

%the level set $S$ is specified as $f_0=f(\bar x)$.

%It is evident that condition \eqref{los} is weaker than
The concept of (first-order) weak sharp minima was first introduced by Ferris in \cite{mcf88}.
%, which accomodates the possibility of non-unique solutions.
 Weak sharp minima are crucial in sensitivity analysis and convergence analysis for various optimization algorithms (see, e.g., \cite{shap,bufe93,bufe95,bumo88,de05,fe91} and the references therein).  It is closely related to the property of linear regularity, metric subregularity and error bound, which has received much attention during the last three decades (see, e.g., \cite{bude02, bude05, bude09,dus10,lmwy11,lhp12, ngzhsiam03,rahmo13, zhngmp09,zhysiam07,zhyjmaa08,zmx12,zhw12} and the references therein). Higher-order weak sharp minima is of great importance in sensitivity analysis of parametric optimization \cite{boni95,shap,dk94}.

Existing research mainly focuses on first-order weak sharp minima, while studies on second-order weak sharp minima are limited. Some relevant works can be found in \cite{shap}.
%\textcolor{blue}{remove or revise the following? Back in the late nineteenth century, Studniarski and Ward \cite{sward99} established sufficient conditions for weak sharp minima of any order $m\geq 1$ for nonconvex functions in finite dimensions through a kind of generalized derivative and generalized tangent cone that involve the level set. In consideration of the importance yet restrictiveness of weak sharp minimizer,} Zheng and Ng proposed the notions of H\"older weak sharp minimizer and H\"older tilt-stable weak minimizer for a proper lower semicontinuous function on a Banach space, and derived optimality conditions for H\"older weak sharp minimizers and stable H\"older weak sharp minimizers and studied their relationships with H\"older metric subregularity/regularity of the subdifferential \cite{zhsiam15,zhnon15}.
{ Let $T_K(\cdot) $
	denote the tangent cone to $K$, $C(x) := \{d \in \mathbb{R}^n | \nabla g(x)d \in T_K(g(x)), \nabla f(x)d \leq 0\}$ denote the
	critical cone at a feasible solution $x$. Let the Lagrange function of problem \eqref{op} be
	$$L(x,\lambda):=f(x)+g(x)^T\lambda.$$ Within the convex analytic framework,  second-order necessary optimality conditions characterizing second-order weak sharp minima are given, as formalized in \cite[Theorem 3.145]{shap}:} Let $\bar x$ be a {local second-order weak sharp minimizer}  for problem \eqref{op} with convex $K$ and suppose that  the Robinson constraint qualification (Robinson CQ) holds at $\bar x$.  Then  for any $\varepsilon\in (0,1/2)$, there exists $\beta >0$ such that for every $x\in S$  sufficiently close to $\bar x$ and $d\in C(x)\cap N_S^{P,\varepsilon}(x)$ where $N_S^{P,\varepsilon}(x)$ denotes the $\varepsilon$-proximal normal cone of $S$ to $ x$, there is a multiplier $\lambda$ fulfilling the first-order necessary optimality condition at $x$ such that
\begin{equation}\label{spc}
\nabla^2_{xx}L( x,\lambda)(d,d)-\sigma_{T_K^2(g( x);\nabla g( x)d)}(\lambda)\geq \beta \|d\|^2,
\end{equation}
where  $\sigma_{C}(\cdot)$ represents the classic support function to set $C$, provided that  the second-order tangent set $T_K^2(g( x);\nabla g( x)d)$ is convex.
%In particular if $K$ and  $T_K^2(g( x);\nabla g( x)d)$ are both convex, then $K(d)$ in  (\ref{spc}) can be taken as $T_K^2(g( x);\nabla g( x)d)$ and (\ref{spc}) becomes the following condition
%\begin{equation}\label{spcnew}
%\nabla^2_{xx}L( x,\lambda)(d,d)-\sigma_{T_K^2(g( x);\nabla g( x)d)}(\lambda)\geq \beta \|d\|^2.
%\end{equation}
Moreover if $K$ is uniformly second-order regular and the property of uniform approximation of critical cones is satisfied at $\bar x$ then the necessary condition (\ref{spc}) becomes sufficient (\cite[Theorem 3.155]{shap}).
% Unfortunately, however,  the second-order tangent set of a convex set is generally not convex and hence even when $K$ is convex, (\ref{spc}) may not hold as a second-order necessary optimality condition.
Even though $K$ is convex, the associated second-order tangent set is often nonconvex, so condition (\ref{spc}) may not necessarily hold as a second-order necessary condition.

% For the case when the constraint set $K$ and  \textcolor{red}{the second-order tangent set} in problem \eqref{op} is convex, there are necessary and sufficient neighborhood conditions (see, e.g., \cite[Theorem 3.145,  \textcolor{red}{ Theorem 3.148 and Theorem 3.155}]{shap}) for the local quadratic growth condition under Robinson's constraint qualification, where the proofs heavily rely on the classical results of convex duality theory (for more details regarding the quadratic growth condition and the nonisolated minima, please refer to the monograph by Bonnans and Shapiro \cite{shap}). A natural question then arises: for the case when the constraint set $K$ in problem \eqref{op} is nonconvex, are there any similar sufficient or necessary optimality conditions for second-order weak sharp minimia?

The primary objective of this paper is to develop new sufficient/necessary optimality conditions for second-order weak sharp minima for the general nonconvex  optimization problem \eqref{op}. Our results do not require the convexity of set $K$ and nonemptyness and/or the convexity of the second-order tangent set $T_K^2(g( x);\nabla g( x)d)$.
In particular in Theorem \ref{cclmby} with $A=\Theta(x,d)$ and $B=\Omega(x;d)$, we have proved the following necessary  conditions for the  second-order weak sharp minima: Let $\bar x$ be a local second-order weak sharp minimizer for problem \eqref{op}.  Then  for any $\varepsilon\in [0,1/2)$, there exists $\beta >0$ such that for every $x\in S$  sufficiently close to $\bar x$ and $d\in C(x)\cap N_S^{P,\varepsilon}(x)$, provided that metric subregularity constraint qualification (MSCQ) holds at $x$ in direction $d$, there exists an M-multiplier $\lambda$ in direction $d$ at $x$ such that $\hat{\sigma}_{T_K^{''}(g( x);\nabla g( x)d)}(\lambda)\leq 0$ and if $T_K^2(g( x);\nabla g( x)d)$ is not empty, then
there is an M-multiplier $\lambda$ in direction $d$ at $x$ such that
\begin{equation}\label{spcnew}
\nabla^2_{xx}L( x,\lambda)(d,d)-\hat{\sigma}_{T_K^2(g( x);\nabla g( x)d)}(\lambda)\geq \beta \|d\|^2,
\end{equation}  where $\hat\sigma_C(\cdot)$  is the so-called lower generalized support function (see Definition \ref{lgsf}), which is generally smaller than the classical support function to set $C$ and coincides with the support function when the set $C$ is convex. And in Corollary \ref{swsm1}, under the directional nondegeneracy condition which is stronger than the directional MSCQ, the lower generalized support function $\hat{\sigma}$ in the above condition can be replaced by the support function $\sigma$. 
Moreover in Corollary \ref{Cor4.5}, we have derived a  sufficient condition which is a point based version of the necessary condition in Theorem \ref{cclmby}. Unlike the classical sufficient condition \cite[Theorem 3.155]{shap}, in our sufficient condition   we do not require extra assumptions such as the uniform second-order regularity of the set $K$ and the property of uniform approximation of the critical cones.

The rest of the paper is organized as follows. Section 2 contains basic notations and preliminary results. In Sections 3, we derive neighborhood necessary optimality conditions and provide examples to illustrate our theoretical results. In Section 4, we develop point-based sufficient second-order optimality conditions. {Conclusions are given in Section 5.}

%%%%%%%%%%%%%%%%%%%%%%%%%%%%%%%%%%%%%%%%%%%%%%%%%%%%%%%%%%%%%%%%%%%%%%%%%%%%%%%%%%%%%%%%%%%%%%%%%%%%%%%%%%%%%%%%%%%%%%%%%%%%%%%%%%%%%%%%%%%%%%%%%%%%%%%%%%%%%%%%%%%%%%%%%%%%%%%%%%%%%%%%%%%

\section{Notations and Preliminary Results}
\label{sec:notation}
In this section, we provide the basic notations and fundamental facts in variational analysis which are used throughout the paper and develop some preliminary results.

Let $S_{\mathbb{R}^n}$ and $B_{\mathbb{R}^n}$ denote the unit sphere and the closed unit ball of $\mathbb{R}^n$, respectively.
$B_r(x)$ represents the closed ball with center $x\in \mathbb{R}^n$ and radius $r>0$. For a set $A\subset \mathbb{R}^n$, we denote its interior, closure, boundary, convex hull and affine hull by ${\rm int}(A), {\rm cl}(A), {\rm bd}(A)$, ${\rm co}(A)$ and ${\rm span}(A)$ respectively. Let $A^\circ$ and $\sigma_A(x)$ stand for the polar cone
and the support function of $A$, respectively, that is, $A^\circ:=\{v\in \mathbb{R}^n \mid \langle v,x\rangle\leq 0\ \text{for\  all}\  x\in A\}$ and $\sigma_A(x):=\sup_{u\in A}\langle x,u\rangle$ for all $x\in \mathbb{R}^n$.
%Let $\dist (x,A):=\inf_{u\in A}\|x-u\|$ denote the point-to-set distance from $x$ to $A$ (in the usual convention, the infimum of the empty set equals $+\infty$). 
{Let $\dist(x, A)$ denote the distance from a point $x$ to the set $A$, defined as $\inf_{u \in A} \|x - u\|$, where the infimum of the empty set is taken to be $+\infty$ by convention.} {Define $P_A(x):=\{w \in {\rm cl}(A) | \|x-w\|=\dist(x,A)\}$ for all $x\in \mathbb{R}^n$.} Let $o: \mathbb{R}_+\rightarrow \mathbb{R}^n$ denote the mapping satisfying $o(t)/t\rightarrow 0$ as $t\downarrow 0$. For $u\in \mathbb{R}^n$, denote by $\{u\}^\bot$ the orthogonal complement of the linear space generated by $u$. For a vector-valued mapping $h:\mathbb{R}^n\rightarrow\mathbb{R}^m$ and sets $B\subset \mathbb{R}^n$, $C\subset \mathbb{R}^m$, we denote $h(B):=\{h(x)|x\in B\}$ and $h^{-1}(C):=\{x| h(x)\in C\}$.

For a vector-valued mapping $g: \mathbb{R}^n\rightarrow \mathbb{R}^m$, we denote by $\nabla g(x)\in \mathbb{R}^{m\times n}$ the Jacobian of $g$ at $x$. The second-order derivative of $g$ at $x$ is denoted by $\nabla^2 g(x)$ and is defined as follows:
$$u^T\nabla^2 g(x):=\lim_{t\rightarrow 0}\frac{\nabla g(x+tu)-\nabla g(x)}{t},\quad\forall u\in \mathbb{R}^n.$$
Hence, we have
$$\nabla^2 g(x)(u,u):=u^T\nabla^2 g(x)u=(u^T\nabla^2 g_1(x)u,\ldots,u^T\nabla^2 g_m(x)u),\quad\forall u\in \mathbb{R}^n.$$
When $g$ is a scalar mapping (i.e. $m=1$), $\nabla^2g(x)$ is identified with the Hessian at $x$.

For a set-valued mapping $M:\mathbb{R}^n\rightrightarrows \mathbb{R}^m$, its graph is defined by
${\rm gph} M:=\{(u,v)\in \mathbb{R}^n\times \mathbb{R}^m \mid v\in M(u)\}$. The inverse mapping $M^{-1}:\mathbb{R}^m\rightrightarrows \mathbb{R}^n$ is defined by $M^{-1}(v)=\{u\in \mathbb{R}^n \mid v\in M(u)\}$ for all $v\in \mathbb{R}^m$.
We denote by $\limsup_{u'\rightarrow u}M(u')$ the
Painlev\'{e}-Kuratowski upper limit, that is,
$$\limsup\limits_{u'\rightarrow u}M(u'):=\{v\in \mathbb{R}^m \mid \exists u_k\rightarrow u,v_k\rightarrow v\,\,{\rm such\, that}\,\,v_k\in M(u_k)\}.$$

For a closed subset $A$ of $\mathbb{R}^n$ and a point $\bar x$ in $A$, the regular/Clarke and Bouligand-Severi tangent/contingent cone to $A$ at $\bar x$ is defined, respectively, by
$$\hat T_A(\bar x):=\liminf\limits_{t\downarrow 0,x\xrightarrow[]{A}\bar x}\frac{A-x}{t}=\{d\in \mathbb{R}^n \mid \forall t_k\downarrow 0, x_k\xrightarrow[]{A}\bar x,\exists d_k\rightarrow d \,\,{\rm with}\,\, x_k+t_kd_k\in A\}$$
and
$$T_A(\bar x):=\limsup\limits_{t\downarrow 0}\frac{A-\bar x}{t}=\{d\in \mathbb{R}^n \mid \exists t_k\downarrow 0, d_k\rightarrow d \,\,{\rm with}\,\, \bar x+t_kd_k\in A\},$$
where $x\xrightarrow[]{A}\bar x$ represents the convergence of $x$ to $\bar x$ with $x\in A$.
The Fr\'{e}chet normal cone of $A$ to $\bar x$ is defined by
$$\hat N_A(\bar x):=\left\{v\in \mathbb{R}^n \mid \limsup_{x\xrightarrow[]{A}\bar x}\frac{\langle v,x-\bar x\rangle}{\|x-\bar x\|}\leq 0\right\}.$$
It is well acknowledged that $\hat N_A(\bar x)=(T_A(\bar x))^\circ$. The proximal normal cone of $A$ to $\bar x$ is defined by
$$N^P_A(\bar x):=\left\{v\in \mathbb{R}^n \mid \exists\tau>0, \,\,s.t.\,\, \bar x+\tau v\in P_A^{-1}(\bar x)\right\}.$$
For a constant $\varepsilon\geq 0$, the $\varepsilon$-proximal normal cone of $A$ to $\bar x$ is defined by
$$N^{P,\varepsilon}_A(\bar x):=\left\{v\in \mathbb{R}^n \mid \dist(v,N^P_A(\bar x))\leq\varepsilon\|v\|\right\}.$$
It is not difficult to verify that $v\in N^P_A(\bar x)$ if and only if there exists $\tau>0$ such that $\bar x+t v\in P_A^{-1}(\bar x)$ for all $t\in (0,\tau)$.
It is clear from the above definitions that the sets $N^P_A(\bar x)$ and $N^{P,\varepsilon}_A(\bar x)$ are cones and $N^{P,0}_A(\bar x)=N^{P}_A(\bar x)$.
Also, since $0\in N^P_A(\bar x)$, we have that $\dist(v,N^P_A(\bar x))\leq \|v\|$,
and hence for $\varepsilon$ greater than $1$, the set $N^{P,\varepsilon}_A(\bar x)$ coincides with the whole space $\mathbb{R}^n$. Therefore,
it makes sense to consider $\varepsilon$-proximal normals for $\varepsilon\in [0,1)$.

Let $N_A(\bar x)$ denote the limiting/Mordukhovich/basic normal cone of $A$ at $\bar x$, that is,
$$N_A(\bar x):=\limsup_{x\xrightarrow[]{A}\bar x}\hat{N}_A(\bar x).$$
{In general, the limiting normal cone lacks convexity, whereas the Fréchet normal cone is always convex by definition (see, for instance, Mordukhovich \cite{boris} and Rockafellar and Wets \cite{rock1}). When the set $A$ is convex, all three types of normal cones—Fréchet, proximal, and limiting—are equivalent to the classical normal cone in convex analysis.

The concept of directional limiting normal cone was first proposed by Ginchev and Mordukhovich \cite{ginc}, and later generalized to arbitrary Banach spaces by Gfrerer in \cite{cfe13}.} For a given direction $d\in \mathbb{R}^n$, the limiting normal cone to $A$ in direction $d$ at $\bar x\in A$ is defined by
\begin{eqnarray*}
N_A(\bar x;d):=\limsup_{t\downarrow 0,d'\rightarrow d}\hat{N}_A(\bar x+td')=\{v\in \mathbb{R}^n \mid  &\exists t_k\downarrow 0, d_k\rightarrow d, v_k\rightarrow v \,\,{\rm with}\\ &v_k\in \hat N_A(\bar x+t_kd_k)\}.
\end{eqnarray*}
It is known that $N_A(\bar x;d)=\emptyset$ if $d\not\in T_A(\bar x)$, $N_A(\bar x;d)\subset N_A(\bar x)$ and $N_A(\bar x;0)= N_A(\bar x)$. Furthermore, if $A$ is convex and $d\in T_A(\bar x)$, we have \cite[Lemma 2.1]{gfer14} that
\begin{equation}\label{ccc}
N_A(\bar x;d)=N_A(\bar x)\cap\{d\}^\bot=N_{T_A(\bar x)}(d).
\end{equation}

Recently, the directional regular/Clarke tangent cone and directional Clarke normal cone were introduced by Gfrerer et al. \cite{ye22}. Let $\hat T_A(\bar x;d)$ denote the directional regular/Clarke tangent cone of $A$ at $\bar x$ in direction $d$, that is,
\begin{equation*}
\aligned
\hat T_A(\bar x;d)&:=\liminf_{t\downarrow 0,d'\rightarrow d, \bar x+td'\in A}T_A(\bar x+td').
\\
\endaligned
\end{equation*}
Let $N^c_A(\bar x;d):={\rm cl}\ {\rm co}\ (N_A(\bar x;d))$ denote the directional Clarke normal cone of $A$ at $\bar x$ in direction $d$.
It is shown in \cite{ye22} that, for a closed set A, $\hat T_A(\bar x;d)$  is closed and convex, $\hat T_A(\bar x)\subset\hat T_A(\bar x;d)$, $\hat T_A(\bar x)=\hat T_A(\bar x;0)$, and furthermore $(\hat T_A(\bar x;d))^\circ=N^c_A(\bar x;d)$.

Next, we review two kinds of second-order tangent sets which play a fundamental role in the second-order analysis later on.

\begin{definition}[Second-order Tangent Sets \cite{shap,pen98}]\label{t1}
Given $A\subset \mathbb{R}^n$, $\bar x\in  A$ and $d \in T_A(\bar x)$.

(i) The outer second-order tangent set to $A$ at $\bar x$ in direction $d$ is defined by
\begin{equation*}
T_A^2(\bar x;d):=\{w\in \mathbb{R}^n \mid \exists t_k\downarrow 0,w_k\rightarrow w\,\,{\rm such\,that}\,\,
\bar x +t_k d+\frac{1}{2}t_k^2w_k\in A\}.
\end{equation*}

(ii) The asymptotic second-order tangent cone to $A$ at $\bar x$ in direction $d$ is defined by
\begin{equation*}
\aligned
T_A^{''}(\bar x;d):=&\{w\in \mathbb{R}^n \mid \exists (t_k,r_k)\downarrow (0,0),w_k\rightarrow w\,\,{\rm such\,that}\\
&t_k/r_k\rightarrow 0,\bar x +t_k d+\frac{1}{2}t_kr_kw_k\in A \}.
\endaligned
\end{equation*}
\end{definition}

{In general, both $T_A^{2}(\bar x;d)$ and $T_A^{''}(\bar x;d)$ are subsets of ${\rm cl}({\rm cone}({\rm cone}(A-\bar x)-d))$. While $T_A^{2}(\bar x;d)$ may not possess a conic structure and can even be empty (see, for instance, Example 3.29 in Bonnans and Shapiro \cite{shap}), $T_A^{''}(\bar x;d)$ is always a cone. Nonetheless, it holds that $$T_A^{2}(\bar x;d) \cup (T_A^{''}(\bar x;d)\backslash \{0\}) \neq \emptyset,\quad\quad \forall d \in T_A(\bar x),$$
which ensures that at least one of these sets is nonempty for any tangent direction $d$; see also \cite[Proposition 2.1]{pen98} and \cite[Proposition 2.2]{OuyangYe}. Especially in the convex setting, we have $T_A^{''}(\bar x;d) = {\rm cl}({\rm cone}({\rm cone}(A-\bar x)-d))$ and $T_A^{2}(\bar x;d)\subset T_A^{''}(\bar x;d)$, implying that the asymptotic second-order tangent cone may strictly contain the second-order tangent set. Additional discussion can be found in \cite{xde} and related literature.} Furthermore, we have the following relationships:

\begin{proposition}\label{ppzz}
Given a closed set $A \subset \mathbb{R}^n$, $\bar x\in A$ and $d\in T_A(\bar x)$. Then,
$$T_A^2(\bar x;d)+\hat T_A(\bar x;d)=T_A^2(\bar x;d)\ \  \text{and}\ \  T_A^{''}(\bar x;d)+\hat T_A(\bar x;d)=T_A^{''}(\bar x;d).$$
\end{proposition}

In what follows, we review the concept of directional neighborhood and some auxiliary results under the directional metric subregularity constraint qualification (directional MSCQ) of the constraint system $g(x)\in K$.

\begin{definition}[Directional Neighborhood\cite{cfe13}]\label{dn}
	Given a direction $d\in \mathbb{R}^n$ and positive numbers $\rho,\delta>0$, the directional neighborhood of direction $d$ is defined as follows:
	\begin{equation*}
		\aligned
		V_{\rho,\delta}(d):&=\left\{w\in \delta B_{\mathbb{R}^n} \mid \left\| \|d\|w-\|w\|d\right\|\leq \rho\|w\|\|d\|\right\}\\
		&=\left\{\begin{matrix}
			\delta B_{\mathbb{R}^n},&{if\,d=0},\\
			\left\{w\in \delta B_{\mathbb{R}^n}\backslash\{0\} \mid \left\|\frac{w}{\|w\|}-\frac{d}{\|d\|}\right\|\leq \rho\right\}\cup \{0\},&{if\,d\not=0}.
		\end{matrix}\right.
		\endaligned
	\end{equation*}
\end{definition}
Note that the directional neighborhood is not a ball except for the case of $d=0$. In general, when $d\in \mathbb{R}^n\backslash\{0\}$, the directional neighborhood is a section of the classical neighborhood.

\begin{definition}[Directional Metric Subregularity {\cite[Definition 1]{cfe13}}  ] \label{subreg} Let  $\bar x$ { be} a feasible point of system $g(x)\in K$.
We say that the metrically subregularity constraint qualification (MSCQ) for system $g(x)\in K$ holds at $\bar x$  in direction $d\in \mathbb{R}^n$, if there are positive numbers $\rho,\delta,\kappa>0$ such that
	\begin{equation*}
		\dist(x,{g^{-1}(K)})\leq \kappa \dist(g(x),K),\quad \forall x\in \bar x+V_{\rho,\delta}(d).
	\end{equation*}
	The infimum of $\kappa$ over all the combinations $\rho,\delta$ and $\kappa$ satisfying the above relation is called the modulus of directional metric subregularity. In the case of $d=0$, we simply say that the metric subregularity constraint
	qualification for {system $g(x)\in K$}  holds at $\bar x$.
\end{definition}
MSCQ is a very weak assumption. In particular according to Robinson \cite[Proposition 1]{robinson}, it is well-known  if $g$ is an affine function and $K$ is a finite union of polyhedral convex sets, then the system $g(x)\in K$ satisfies MSCQ at any feasible point.  It is also well-known that Robinson's CQ is a sufficient condition for MSCQ. The following 
sufficient conditions {for} MSCQ which are weaker than Robinson's CQ are given in \cite[Corollary 1]{gferk16}. 

\begin{proposition} [cf.{\cite[Corollary 1]{gferk16}}]
	Given {a system $g(x)\in K$} where $g:\mathbb{R}^n\rightarrow \mathbb{R}^m$ is continuously differentiable and $K$ is a closed set.  Then,  MSCQ for {system $g(x)\in K$} holds at a feasible point $\bar x$ {in direction $d$} if one of the following conditions is fulfilled:
	\begin{itemize}
		\item[(i)] First-order sufficient condition for metric subregularity (FOSCMS) {at $\bar x$ in direction $0 \neq d \in \mathbb{R}^n$ with $\nabla g(\bar x)d \in T_{K}(g(\bar x))$ }: one has
		$$\nabla g(\bar x)^T\lambda = 0, \lambda \in N_K(g(\bar x);\nabla g(\bar x)d) \Rightarrow  \lambda = 0.$$
		\item[(ii)] Second-order sufficient condition for metric subregularity (SOSCMS) {at $\bar x$ in direction $0 \neq d \in \mathbb{R}^n$ with $\nabla g(\bar x)d \in T_{K}(g(\bar x))$}: $g$ is twice Fréchet differentiable at $\bar x$, $K$ is the union of finitely many convex polyhedra, and 
			$$\nabla g(\bar x)^T\lambda = 0, \lambda \in N_K(g(\bar x);\nabla g(\bar x)d), d^T \nabla^2(\lambda^Tg)(\bar x)d \geq 0 \Rightarrow  \lambda = 0.$$
		\end{itemize}
\end{proposition}

Other sufficient conditions for MSCQ can be found  in e.g., \cite{yezhou18}.

The following result is important in the proof of our main result in {Section 3}.
\begin{proposition}[cf.{\cite[Proposition 5]{ye22},\cite[Proposition 2.7]{OuyangYe}}]\label{2.2}
	Given $\bar x\in \Phi:=g^{-1}(K)$ and $d\in T_\Phi(\bar x)$. Then
	\begin{equation}\label{2cb}
		T_\Phi^2(\bar x;d)\subset \{w\in \mathbb{R}^n \mid \nabla g(\bar x)w+\nabla^2 g(\bar x)(d,d)\in T_K^2(g(\bar x);\nabla g(\bar x)d)\}
	\end{equation}
	and
	\begin{equation}\label{2cb1}
		T_\Phi^{''}(\bar x;d)\subset \{w\in \mathbb{R}^n \mid \nabla g(\bar x)w\in T_K^{''}(g(\bar x);\nabla g(\bar x)d)\}.
	\end{equation}
	If, in addition, assume that MSCQ holds at $\bar x$ in direction $d\in \mathbb{R}^n$ for the constraint system $g(x)\in K$ with modulus $\kappa$, then inclusions \eqref{2cb} and \eqref{2cb1} hold as equality.
\end{proposition}

%%%%%%%%%%%%%%%%%%%%%%%%%%%%%%%%%%%%%%%%%%%%%%%%%%%%%%%%%%%%%%%%%%%%%%%%%%%%%%%%%%%%%%%%%%%%%%%%%%%%%%%%%%%%%%%%%%%%%%%%%%%%%%%%%%%%%%%%%%%%%%%%%%%%%%%%%%%%%%%%%%%%%%%%%%%%%%%%%%%%%%%%%%%

\section{Necessary optimality conditions for second-order weak sharp minima}
\label{sec:sows}
In this section, we focus on developing neighborhood necessary optimality conditions for second-order weak sharp minima of general nonconvex optimization problem \eqref{op}, in both implicit and explicit form. 

In \cite[Theorem 4.2]{OuyangYe}, the authors established necessary optimality conditions for local minima. { Notably, second-order weak sharp minima constitute a special subclass of local minima. } Motivated by their result, we next present necessary optimality conditions for second-order weak sharp minima. To this end, we first introduce the following lemma.

\begin{lemma} [cf.{ \cite[Lemma 3.144]{shap} }]\label{shap3.144}
For $\epsilon \geq 0$ let $d$ be a $\epsilon$-proximal normal to $S$ at $x \in S$. Then for all $t>0$ small enough,
$$\dist(x+td,S) \geq t(1-2\epsilon)\|d\|.$$
\end{lemma}
By exploiting the properties of outer second-order tangent set and asymptotic second-order tangent cone to the feasible set $\Phi$ and utilizing the above estimates of the distance to $S$ along $\varepsilon$-proximal normals, we derive the following implicit necessary condition for second-order weak sharp minima.

\begin{theorem}\label{nscll}
Let $\bar x$ be a local second-order weak sharp minimizer of problem \eqref{op} with the corresponding constants $\kappa >0$ and $\delta >0$. Then for any $\varepsilon\in [0,1/2)$, for every $x\in S\cap B_{\delta}(\bar x)$, $d\in C(x)\cap N_S^{P,\varepsilon}(x)$, and any $\lambda \in \mathbb{R}^m$ with $\nabla_x L( x, \lambda)=0$, one has
\begin{enumerate}[{\rm (i)}]
\item $\sigma_{\nabla g(x)(T_\Phi^{''}(x,d))}(\lambda)\leq 0$;
\item $\nabla_{xx}^2L(x,\lambda)(d,d)-\sigma_{\nabla g(x)(T_\Phi^2(x;d))+\nabla^2 g(x)(d,d)}(\lambda)\geq 2\kappa(1-2\varepsilon)^2\|d\|^2$.
\end{enumerate}
\end{theorem}

\beginproof
Pick any $\varepsilon\in [0,1/2)$, $x\in S\cap B_{\delta}(\bar x)$, $d\in C(x)\cap N_S^{P,\varepsilon}(x)$, and $\lambda \in \mathbb{R}^m$ with $\nabla_x L( x, \lambda)=0$ and let them be fixed. Since $x\in S$ is automatically a local minimizer of problem \eqref{op}, we deduce from \cite[Theorem 4.2]{OuyangYe} that assertion (i) holds true.
% and
%\begin{equation}\label{xxzq}
%\nabla_{xx}^2L( x,\lambda)(d,d)-\sigma_{\nabla g(x)(T_\Phi^2(x;d))+\nabla^2 g(x)(d,d)}(\lambda)= \alpha\geq0,
%\end{equation} where
%$$\alpha:=\inf_{w\in T_{\Phi}^2( x;d)}(\nabla f( x)w+\nabla^2 f( x)(d,d)).$$

Now pick arbitrary $w\in T_\Phi^2(x;d)$ and obtain sequences $t_k\downarrow 0$ and $w_k\rightarrow w$ such that
$\tilde{x}_k:=x +t_k d+\frac{1}{2}t_k^2w_k\in \Phi$ for all $k\in \mathbb{N}$. It then follows from $d\in C(x)$ and Taylor expansion that
\begin{equation}\label{te3}
f(\tilde{x}_k)-f(x)\leq \frac{1}{2}t_k^2(\nabla f(x)w_k+\nabla^2 f(x)(d,d))+o(t_k^2).
\end{equation}
Since $d\in N_S^{P,\varepsilon}(x)$, we obtain from Lemma \ref{shap3.144} that, for sufficiently large $k$,
\begin{equation*}
\aligned
\dist (\tilde{x}_k, S)&=\dist (x +t_k d+\frac{1}{2}t_k^2w_k, S)\\
&\geq \dist (x +t_k d, S)-\frac{1}{2}t_k^2\|w_k\|\\
&\geq t_k(1-2\varepsilon)\|d\|-\frac{1}{2}t_k^2\|w_k\|.
\endaligned
\end{equation*}
Therefore, for sufficiently large $k$, we have
\begin{equation}\label{te4}
f(\tilde{x}_k)-f(x)=f(\tilde{x}_k)-f(\bar x)\geq \kappa \dist^2 (\tilde{x}_k, S) \geq \kappa t_k^2(1-2\varepsilon)^2\|d\|^2+o(t_k^2).
\end{equation}
Combining inequalities \eqref{te3} and \eqref{te4} and passing to the limit, we obtain that
$$\nabla f(x)w+\nabla^2 f(x)(d,d)\geq 2\kappa(1-2\varepsilon)^2\|d\|^2.$$
Then for any $\lambda \in \mathbb{R}^m$ with $\nabla_x L( x, \lambda)=0$, we have
\begin{equation*}
\aligned
&\sigma_{\nabla g(x)(T_\Phi^2(x;d))+\nabla^2 g(x)(d,d)}(\lambda)=\sup_{u\in \nabla g(x)(T_\Phi^2(x;d))+\nabla^2 g(x)(d,d)}\langle\lambda,u\rangle\\
&=\sup_{w\in T_\Phi^2(x;d)}\langle\lambda,\nabla g(x)w+\nabla^2 g(x)(d,d)\rangle
\\&=\sup_{w\in T_\Phi^2(x;d)}\langle\nabla g(x)^T\lambda,w\rangle+\langle\lambda,\nabla^2 g(x)(d,d)\rangle
\\&=\sup_{w\in T_\Phi^2(x;d)}-(\nabla f(x)w+\nabla^2 f(x)(d,d))+\nabla^2 f(x)(d,d)+\langle\lambda,\nabla^2 g(x)(d,d)\rangle
\\&=-\inf_{w\in T_\Phi^2(x;d)}(\nabla f(x)w+\nabla^2 f(x)(d,d))+\nabla^2_{xx} L(x,\lambda)(d,d)
\\&\leq -2\kappa(1-2\varepsilon)^2\|d\|^2+\nabla^2_{xx} L(x,\lambda)(d,d),
\endaligned
\end{equation*}
which indicates that assertion (ii) holds true.
\endproof

{Comparing Theorem \ref{nscll} with \cite[Theorem 4.2]{OuyangYe}, the first part (i) is identical in both, whereas the second part (ii) in Theorem \ref{nscll} is more stringent due to the requirement that  $\kappa>0$.}

 In the following result, we provide a neighborhood necessary optimality condition for second-order weak sharp minima of problem \eqref{op} that are valid for all critical directions which are not necessarily $\varepsilon$-proximal normals.

\begin{theorem}\label{nsc22}
Let $\bar x$ be a local second-order weak sharp minimizer of problem \eqref{op} with the corresponding constants $\kappa >0$ and $\delta >0$. Then for every $x\in S\cap B_{\delta}(\bar x)$, $d\in C(x)$, and any $\lambda \in \mathbb{R}^m$ with $\nabla_x L( x, \lambda)=0$, one has
\begin{enumerate}[{\rm (i)}]
\item $\sigma_{\nabla g(x)(T_\Phi^{''}(x,d))}(\lambda)\leq 0$;
\item $\nabla_{xx}^2L(x,\lambda)(d,d)-\sigma_{\nabla g(x)(T_\Phi^2(x;d))+\nabla^2 g(x)(d,d)}(\lambda)\geq 2\kappa \dist(d,T_S(x))^2$.
\end{enumerate}
\end{theorem}

\beginproof
Pick any $x\in S\cap B_{\delta}(\bar x)$, $d\in C(x)$, and any $\lambda \in \mathbb{R}^m$ with $\nabla_x L( x, \lambda)=0$. It follows from \cite[Theorem 4.2]{OuyangYe} that assertion (i) holds.

To establish (ii), we pick arbitrary $w\in T_\Phi^2(x;d)$ and obtain sequences $t_k\downarrow 0$ and $w_k\rightarrow w$ such that
$\tilde{x}_k:=x +t_k d+\frac{1}{2}t_k^2w_k\in \Phi$ for all $k\in \mathbb{N}$. By \cite[Lemma 3.147]{shap}, we have
\begin{equation*}
\aligned
\dist (\tilde{x}_k, S)&=\dist (x +t_k d+\frac{1}{2}t_k^2w_k, S)\\
&\geq t_k \dist (d+\frac{1}{2}t_kw_k, T_S(x))+ o(t_k)\\
&\geq t_k\dist (d, T_S(x))-\frac{1}{2}t_k^2\|w_k\|+ o(t_k).
\endaligned
\end{equation*}
Then we deduce that
\begin{equation}\label{te5}
f(\tilde{x}_k)-f(x)=f(\tilde{x}_k)-f(\bar x)\geq \kappa \dist^2 (\tilde{x}_k, S) \geq \kappa t_k^2\dist^2(d, T_S(x))+o(t_k^2)
\end{equation}
holds for sufficiently large $k$. Combining \eqref{te3} with \eqref{te5} and passing to the limit, we obtain
$$\nabla f(x)w+\nabla^2 f(x)(d,d)\geq 2\kappa \dist(d,T_S(x))^2.$$
Then, in analogous to the proof of Theorem \ref{nscll} (ii), we arrive at assertion (ii).
\endproof

Compared to Theorem \ref{nscll}, Theorem \ref{nsc22} allows a broader range of directions $d$ (not limited to $\epsilon$-proximal normals), while the right-hand side is replaced by a distance function. Theorem \ref{nsc22} recovers \cite[Theorem 4.2]{OuyangYe} when $\bar x$
is an isolated local minimizer.

{We present the following examples to illustrate the optimality conditions derived above. In the first example, we demonstrate how to apply these conditions to a specific problem. In the second example, we show how the necessary optimality conditions can be used to exclude points that are not second-order weak sharp minima. }

\begin{example}
Consider the optimization problem
\begin{equation*}
\aligned
\min &\ f(x_1,x_2):=x_2^2\\
{\rm s.t.} &\ g(x_1,x_2):=x_1^2-2x_1+x_2^2\in K:= [-3/4,0{\color{blue}]}
\endaligned
\end{equation*}
at $\bar x=(0,0)\in \mathbb{R}^2$. Then it is easy to calculate that the feasible set
 $\Phi=g^{-1}(K)=\{(x_1,x_2)\in \mathbb{R}^2 \mid 1/4\leq(x_1-1)^2+x_2^2\leq 1\}.$
Consider the set $S=[0,1/2]\times \{0\}\subset \mathbb{R}^2$. Observe that $T_S(\bar x)=\mathbb{R}_+\times \{0\}$ and $f$ is constant on $S$. Furthermore, we can see that $\bar x$ is a second-order weak sharp minimizer of the above nonconvex optimization problem.

%(i) Observe that $T_S(\bar x)=\mathbb{R}_+\times \{0\}$. Let $d=(1,0)\in T_S(\bar x)$, we have $T_S^{''}(\bar x;d)=T_S^{2}(\bar x;d)=\mathbb{R}\times \{0\}$.
%And then, $\nabla f(\bar x)v=\nabla f(\bar x)v+\nabla^2 f(\bar x)(d,d)=0$ for any $v\in T_S{''}(\bar x;d)$. Now pick any $x=(x_1,x_2) \in S\backslash \{0\}$ and $d\in T_S(x)$,  it easy to verify that assertions (i)-(iii) of Proposition \ref{nsop} hold true.

(i) In order to illustrate Theorem \ref{nscll}, for convenience, we consider $\varepsilon=0$ and $\delta=1/4$.
Let $d=(0,1)\in C(\bar x)\cap N_S^P(\bar x)=\{0\}\times \mathbb{R}$,
we have $T_\Phi^{''}(\bar x;d)=\{(x_1,x_2)\in \mathbb{R}^2 \mid x_1\geq 0\}$ and $T_\Phi^{2}(\bar x;d)=\{(x_1,x_2)\in \mathbb{R}^2 \mid x_1\geq 1\}$.
{Solving  $\nabla_x L( \bar x, \lambda)=0$, we obtain $\lambda = 0$. Then,}
$\sigma_{\nabla g(\bar x)(T_\Phi^{''}(\bar x,d))}(\lambda) = 0$ and
$\nabla_{xx}^2L(\bar x,\lambda)(d,d)-\sigma_{\nabla g(\bar x)(T_\Phi^2(\bar x;d))+\nabla^2 g(\bar x)(d,d)}(\lambda)=\nabla^2f(\bar x)(d,d)=2= 2\|d\|^2$.
Therefore, assertions (i) and (ii) of Theorem \ref{nscll} hold for $x=\bar x$ and $d=(0,1)$.
For other choices of $x\in S$ and $d\in C(x)\cap N_S^P(x)$, we can verify the assertions analogously.

(ii) Let $d=(1,1)\in C(\bar x)=\mathbb{R}_+\times \mathbb{R}$, then $\dist(d,T_S(\bar x))=1$.
{Solving  $\nabla_x L( \bar x, \lambda)=0$, we obtain $\lambda = 0$. Then,} $\nabla_{xx}^2L(\bar x,\lambda)(d,d)-\sigma_{\nabla g(\bar x)(T_\Phi^2(\bar x;d))+\nabla^2 g(\bar x)(d,d)}(\lambda)=\nabla^2 f(\bar x)(d,d)= 2 = 2\dist(d,T_S(\bar x))^2$. Therefore, conclusions (i) and (ii) of Theorem \ref{nsc22} hold for $x=\bar x$ and $d=(1,1)$. For other choices of $x\in S$ and $d\in C(x)$, it can be verified in a similar manner.

\end{example}

Next, we use second-order necessary optimality conditions in Theorem \ref{nscll} (or Theorem \ref{nsc22}) to rule out points that are not second-order weak sharp minima.
\begin{example}[cf.{\cite[Example 2]{ye22} }]
	Consider the optimization problem
	\begin{equation*}
		\aligned
		\min_{x \in \mathbb{R}} &\ f(x):=-\frac{1}{2} x^2\\
		{\rm s.t.} &\ g(x):=(x^2,x)\in K:= C_1 \cup C_2
		\endaligned
	\end{equation*}
	at $\bar x=0$, where $C_1:=\{(x_1,x_2)\in \mathbb{R}^2 \mid (x_1-1)^2+x_2^2\leq 1\}$ and $C_2:=\{(x_1,x_2)\in \mathbb{R}^2 \mid (x_1+1)^2+x_2^2\leq 1\}$. Then it is easy to calculate that the feasible set
	$\Phi=g^{-1}(K)=[-1,1]$ and thus $\bar x$ is not a second-order weak sharp minimizer. We {now} check that our results can reject $\bar x$ as a second-order weak sharp minimizer.
	Consider the set $S=\{0\}$, $\varepsilon=0$ and $\delta=1/4$.
	Let $d=1\in C(\bar x)\cap N_S^P(\bar x)= \mathbb{R}$,
	we have  $T_\Phi^{2}(\bar x;d)=\mathbb{R}$.
	There exists $\lambda=(0,0)\in N_K(g(\bar x))$ such that
	$\sigma_{\nabla g(\bar x)(T_\Phi^{''}(\bar x,d))}(\lambda) = 0$ and
	$\nabla_{xx}^2L(\bar x,\lambda)(d,d)-\sigma_{\nabla g(\bar x)(T_\Phi^2(\bar x;d))+\nabla^2 g(\bar x)(d,d)}(\lambda)=\nabla^2f(\bar x)(d,d)=-1<0$.
	Therefore, by Theorem \ref{nscll} (or Theorem \ref{nsc22}), $\bar x$ is not a second-order weak sharp minimizer.
	\end{example}

Recall that in the recent work \cite{ye22}, the authors introduced the concept of lower generalized support function as follows:
\begin{definition}[cf.{\cite[Definition 7]{ye22}}]\label{lgsf}
Let $S \subset \mathbb{R}^n$ be a nonenmpty closed set. For every subset $A \subset \mathbb{R}^n$, the lower generalized support function to $S$ with respect to $A$ is defined as the mapping $\hat{\sigma}_{S, A}: \mathbb{R}^n \rightarrow \mathbb{R} \cup\{ \pm \infty\}$ by
$$
\hat{\sigma}_{S, A}(\lambda):=\liminf _{\lambda^{\prime} \rightarrow \lambda} \inf _u\left\{\left\langle\lambda^{\prime}, u\right\rangle \mid u \in N_S^{-1}\left(\lambda^{\prime}\right) \cap A\right\} .
$$
If $S=\emptyset$, then we define $\hat{\sigma}_{S, A}(\lambda)=-\infty$ for all $\lambda$. When $A=\mathbb{R}^n$, we use $\hat{\sigma}_S$ in place of $\hat{\sigma}_{S, \mathbb{R}^n}$.
\end{definition}

Note that by the definition, we have $\hat{\sigma}_S \leq \hat{\sigma}_{S, A}$ for every subset $A \subset \mathbb{R}^n$ and $\hat{\sigma}_{S, B} \leq \hat{\sigma}_{S, A}$ whenever $A \subset B \subset \mathbb{R}^n$. We also observe that $\hat{\sigma}_{S, A}(\lambda)=+\infty$ if for all $\lambda^{\prime}$ close to $\lambda, N_S^{-1}\left(\lambda^{\prime}\right) \cap A=\emptyset$. The lower generalized support function is always less than or equal to the support function and that both functions coincide when the underlying set is convex \cite[Proposition 6]{ye22}.

Note that the necessary optimality conditions in {Theorems 3.2 and 3.3} are formulated by using the feasible solutions and hence is implicit. Next, we develop explicit form of neighborhood second-order necessary optimality conditions for second-order weak sharp minima in terms of directional M-multipliers under MSCQ. Recall that the directional Mordukhovich (M-) multiplier set is defined as:
$$\Lambda (\bar x, d):= \{\lambda \in \mathbb{R}^m \mid \nabla_xL(\bar x, \lambda)=0, \lambda \in N_K(g(\bar x);\nabla g(\bar x)d)\}.$$
Given $x\in S$ and $d \in T_\Phi(x)$. Define
$$\Theta(x;d):=\nabla g(x)(T_\Phi^{''}(x,d)), \quad \Omega(x;d):=\nabla g(x)(T_\Phi^2(x;d))+\nabla^2 g(x)(d,d).$$
Observe that, if MSCQ holds at $\bar x\in S$ for the constraint system $g(x)\in K$, then there exists $\delta'>0$ such that MSCQ holds at all $x\in S\cap B_{\delta'}(\bar x)$.
Therefore, we can apply \cite[Proposition 4.6]{OuyangYe} with $\Theta(x;d)$ and $\Omega(x;d)$. {Note that a point $\bar x \in \Phi$ is said to be a local optimal solution of problem \eqref{op} in direction $d \in \mathbb{R}^n$, if there exist positive numbers $\rho, \delta>0$ such that $f(x) \geq f(\bar x)$ for all $x \in \Phi \cap (\bar x+V_{\rho, \delta}(d))$ \cite[Definition 3.1]{OuyangYe}.}

\begin{lemma}[cf.{\cite[Proposition 4.6]{OuyangYe}}] Let $d \in T_{\Phi}(x)$ and $x$ be a local optimal solution of problem \eqref{op} in direction d. Suppose that $\nabla f(x) d=0$ and $M S C Q$ holds in direction $d$ for the constraint system $g(x) \in K$. Then the following relationships hold:
	\begin{itemize}
\item[(i)] There exists $\lambda \in \Lambda(x ; d)$ such that for any $A \supset \Theta(x;d)$, we have

$$\hat\sigma_{T^{''}_K (g( x);\nabla g( x)d)}(\lambda) \leq \hat\sigma_{T^{''}_K (g( x);\nabla g( x)d),A}(\lambda) \leq \sigma_{\Theta(x;d)}(\lambda),$$

\item[(ii)] If $T_K^2(g(x) ; \nabla g(x) d) \neq \emptyset$, then there exists $\lambda \in \Lambda(x ; d)$ such that, for any $B \supset \Omega(x;d)$, we have

$$\hat\sigma_{T^2_K (g( x);\nabla g( x)d)}(\lambda) \leq \hat\sigma_{T^2_K (g( x);\nabla g( x)d),B}(\lambda) \leq \sigma_{\Omega(x;d)}(\lambda).$$
\end{itemize}
\end{lemma}

We can deduce from Theorem \ref{nscll} and Theorem \ref{nsc22}, respectively, that the following explicit necessary optimality conditions hold. {Note that, under the assumption of MSCQ, every direction $d \in C(x)$ satisfies $\nabla f(x) d=0$.}
	
	\begin{theorem}\label{cclmby}
 Let $\bar x$ be a local second-order weak sharp minimizer of problem \eqref{op} with the corresponding constants $\kappa >0$ and $\delta >0$. Suppose that  MSCQ holds at $\bar x$. Then for any $\varepsilon\in [0,1/2)$, there exists $\delta'\in (0,\delta)$ such that for every $x\in S\cap B_{\delta'}(\bar x)$, $d\in C(x)\cap N_S^{P,\varepsilon}(x)$, the following assertions hold:
\begin{enumerate}[{\rm (i)}]
\item There exists $\lambda\in \Lambda( x;d)$ such that for any set $A\supset \Theta(x;d)$, we have $\hat\sigma_{T^{''}_K (g( x);\nabla g( x)d),A}(\lambda)\leq  0$. In particular,
$$\hat\sigma_{T^{''}_K (g( x);\nabla g( x)d)}(\lambda)\leq  0;$$
\item If $T^2_K (g( x);\nabla g( x)d)\not=\emptyset$, then there exists $\lambda\in \Lambda( x;d)$ such that, for any set $B\supset \Omega(x;d)$, we have
\begin{equation}\label{aayjj}
\nabla_{xx}^2L(x,\lambda)(d,d)-\hat\sigma_{T^2_K (g( x);\nabla g( x)d),B}(\lambda)\geq 2\kappa(1-2\varepsilon)^2\|d\|^2.
\end{equation}
%where $\alpha$ is defined as in Theorem \ref{nss}.
In particular,
\begin{equation}\label{aayjjnew}
\nabla_{xx}^2L( x,\lambda)(d,d)-\hat\sigma_{T^2_K (g( x);\nabla g( x)d)}(\lambda)\geq  2\kappa(1-2\varepsilon)^2\|d\|^2.
\end{equation}
\end{enumerate}
Besides, for any $x\in S\cap B_{\delta'}(\bar x)$ and $d\in C(x)$, we also have (i) and (ii) hold true with \eqref{aayjj} and \eqref{aayjjnew} being replaced by
\begin{equation}\label{}
	\nabla_{xx}^2L(x,\lambda)(d,d)-\hat\sigma_{T^2_K (g( x);\nabla g( x)d),B}(\lambda)\geq 2\kappa \dist(d,T_S(x))^2
\end{equation}
and
\begin{equation}\label{}
	\nabla_{xx}^2L( x,\lambda)(d,d)-\hat\sigma_{T^2_K (g( x);\nabla g( x)d)}(\lambda)\geq  2\kappa \dist(d,T_S(x))^2,
	\end{equation} respectively.
\end{theorem}

{In the above, we derived {explicit} optimality conditions using the lower generalized support function, which is {weaker} than the standard case. Next, we aim to establish explicit optimality conditions using the classical support function. Motivated by \cite[Corollary 4.10]{OuyangYe}, we first introduce the following multipliers.} Consider the following directional Clarke (C-) multiplier set:
$$\Lambda^c (\bar x; d):= \{\lambda \mid \nabla_xL(\bar x, \lambda)=0, \lambda \in N^c_\Lambda(g(\bar x);\nabla g(\bar x)d)\}$$
and the directional Robinson's constraint qualification (DirRCQ) in direction $d$:
\begin{equation}\label{dircq}
\nabla g(\bar x)^T\lambda=0,\lambda\in N^c_K(g(\bar x);\nabla g(\bar x)d)\Rightarrow \lambda=0.
\end{equation}
{Clearly, the directional C-multiplier set is a closed and convex set, and it typically contains the directional M-multiplier set. As a result, according to \cite[Theorem 1]{gferk16}, the DirRCQ condition \eqref{dircq} imposes a stronger requirement than the MSCQ in the direction $d$. It is also important to emphasize that both the directional M- and C-multiplier sets, as well as the DirRCQ and the directional MSCQ, represent weaker versions of their corresponding nondirectional counterparts, which can be regarded as the particular case when $d=0$.}

In what follows, we derive a necessary condition for second-order weak sharp minima of problem \eqref{op} in terms of directional C-multipliers under DirRCQ.

\begin{theorem}\label{swmc}
Let $\bar x$ be a local second-order weak sharp minimizer of problem \eqref{op} with the corresponding constants $\kappa >0$ and $\delta >0$. Suppose that DirRCQ holds at every $x \in B_\delta (\bar x)$ in direction $d\in C(x)$. Then for any $\varepsilon\in [0,1/2)$, $x\in S\cap B_{\delta}(\bar x)$, $d\in C(x)\cap N_S^{P,\varepsilon}(x)$, the following assertions hold:
\begin{enumerate}[{\rm (i)}]
\item For every $v \in T_K^{''}(g(x);\nabla g(x)d)$, there exists $\lambda_v \in \Lambda^c (x; d)$ such that $\langle\lambda_v,v\rangle\leq 0$;
\item For every $w \in T_K^2(g(x);\nabla g(x)d)$, there exists $\lambda_w \in \Lambda^c (x; d)$ such that
\begin{equation}\label{nonprox}
\nabla_{xx}^2L(x,\lambda_w)(d,d)-\langle\lambda_w,w\rangle\geq 2\kappa(1-2\varepsilon)^2\|d\|^2.
\end{equation}
\end{enumerate}

Besides, for any $x\in S\cap B_{\delta'}(\bar x)$ and $d\in C(x)$, we also have (i) and (ii) hold true with \eqref{nonprox} being replaced by
\begin{equation}\label{nonprox1}
\nabla_{xx}^2L(x,\lambda_w)(d,d)-\langle\lambda_w,w\rangle\geq 2\kappa\dist(d,T_S(x))^2.\end{equation}
\end{theorem}

\beginproof
Pick any $\varepsilon\in [0,1/2)$, $x\in S\cap B_{\delta}(\bar x)$, $d\in C(x)\cap N_S^{P,\varepsilon}(x)$ and let them be fixed. Since DirRCQ imiplies MSCQ for the system $g(x)\in K$, by { Proposition \ref{2.2}}, we have
$$T_\Phi^{''}(x;d)= \{u\in \mathbb{R}^n \mid \nabla g(x)u\in T_K^{''}(g(x);\nabla g(x)d)\}$$
and
$$T_\Phi^{2}(x;d)= \{w\in \mathbb{R}^n \mid \nabla g(x)w+\nabla^2 g(x)(d,d)\in T_K^{''}(g(x);\nabla g(x)d)\}.$$
The rest of the proof takes the idea in the arguments of Gfrerer et al.\cite[Proposition 8]{ye22}.

For (i), we take $v\in T_K^{''}(g(x);\nabla g(x)d)$. Then it follows from Proposition \ref{ppzz} and the proof of Theorem \ref{nscll} (i) that, the following conic linear program
\begin{equation}\label{clp}
\min_{u} \nabla f(x)u\,\,\, {\rm s.t.}\,\,\,\nabla g(x)u\in v+\hat T_K(g(x);\nabla g(x)d)
\end{equation}
has nonnegative optimal value and $\Lambda^c(x;d)$ is nonempty. Since $(\hat T_K(g(x);\nabla g(x)d))^\circ=N^c_K(g(x);\nabla g(x)d)$,  the dual program of \eqref{clp} can be written as
\begin{equation}\label{clp1}
\max_{\lambda\in \Lambda^c(x;d)} -\lambda^Tv.
\end{equation}
Note that $\Lambda^c(x;d)$ is compact and the DirRCQ implies that
$$0\in {\rm int}(\nabla g(x)\mathbb{R}^n-v-\hat T_K(g(x);\nabla g(x)d))$$
thanks to \cite[Lemma 6]{ye22},
we deduce from Bonnans and Shapiro \cite[Theorem 2.187]{shap} that, there is no dual
gap between problems \eqref{clp} and \eqref{clp1}, and the dual program \eqref{clp1} has an optimal solution $\lambda_v$ such that
$$-\langle\lambda_v,v\rangle=\max_{\lambda\in \Lambda^c(x;d)} -\lambda^Tv\geq 0.$$
This shows (i) holds.

To establish (ii), take $w\in T_K^{2}(g(x);\nabla g(x)d)$. Then we deduce from the proof of Theorem \ref{nscll} (ii) that the following conic linear program
\begin{equation}\label{clp2}
\aligned
\min_{u}\quad &\nabla f(x)u+\nabla^2 f(x)(d,d)-2\kappa(1-2\varepsilon)^2\|d\|^2\\
{\rm s.t.}\quad &\nabla g(x)u+\nabla^2 g(x)(d,d)\in w+\hat T_K(g(x);\nabla g(x)d)
\endaligned
\end{equation}
has nonnegative optimal value. The dual program of the conic linear program \eqref{clp2} is
\begin{equation}\label{dclp2}
\aligned
\max_{\lambda\in \Lambda^c(x;d)}\quad &\nabla_{xx}^2L(x,\lambda)(d,d)-\lambda^Tw-2\kappa(1-2\varepsilon)^2\|d\|^2\\
{\rm s.t.}\quad &\nabla g(x)^T\lambda+\nabla f(x)=0.
\endaligned
\end{equation}
It is easy to see from \cite[Lemma 6]{ye22} that the DirRCQ implies
$$0\in {\rm int}(\nabla g(x)\mathbb{R}^n-\nabla^2 g(x)(d,d)-w-\hat T_K(g(x);\nabla g(x)d{\color{blue})).}$$
Consequently, we deduce from Bonnans and Shapiro \cite[Theorem 2.187]{shap} that, there is no dual
gap between problems \eqref{clp2} and \eqref{dclp2}, and the dual program \eqref{dclp2} has an optimal solution $\lambda_w$ such that the inequality in assertion (ii) holds.

For critical directions which are not necessarily $\varepsilon$-proximal normals, by applying the results of Theorem \ref{nsc22} instead of Theorem \ref{nscll}, we can also derive (i) and (ii) with \eqref{nonprox1} being true.
\endproof

In the following corollary, we show that under directional Robinson constraint qualification which is weaker than Robinson CQ, we obtain stronger necessary optimality condition in comparison with {\cite[Theorems 3.145 and 3.148]{shap}} in that, the sharper directional Clarke Multipliers are employed, the second-order asymptotic tangent cone is considered in the case of $T_K^2(g(x);\nabla g(x)d) \neq \emptyset$, and most importantly, the constraint set $K$ and the second-order tangent set $T_K^2(g(x);\nabla g(x)d)$ here can be nonconvex.

\begin{corollary}\label{swsm}
Let $\bar x$ be a local second-order weak sharp minimizer of problem \eqref{op} with the corresponding constants $\kappa >0$ and $\delta >0$. Suppose that DirRCQ holds at every $x \in B_\delta (\bar x)$ in direction $d\in C(x)$. Then for any $\varepsilon\in [0,1/2)$, $x\in S\cap B_{\delta}(\bar x)$, $d\in C(x)\cap N_S^{P,\varepsilon}(x)$, the following assertions hold:
\begin{enumerate}[{\rm (i)}]
\item For any nonempty convex subset $K(d)\subset T^{''}_K (g(x);\nabla g(x)d)$, there exists $\lambda \in \Lambda^c (x; d)$ such that $\sigma_{K(d)}(\lambda)\leq 0$;
\item If $T_K^2(g(x);\nabla g(x)d) \neq \emptyset$, then for any nonempty convex subset $\hat K(d)\subset T^2_K (g(x);\nabla g(x)d)$, there exists $\lambda \in \Lambda^c (x; d)$ such that
\begin{equation}\label{prox2}
\nabla_{xx}^2L(x,\lambda)(d,d)-\sigma_{\hat K(d)}(\lambda)\geq 2\kappa(1-2\varepsilon)^2\|d\|^2.
\end{equation}
\end{enumerate}

Besides, for any $x\in S\cap B_{\delta}(\bar x)$ and $d\in C(x)$, we also have (i) and (ii) hold true with \eqref{prox2} being replaced by
\begin{equation}\label{nonprox2}
\nabla_{xx}^2L(x,\lambda)(d,d)-\sigma_{\hat K(d)}(\lambda)\geq 2\kappa\dist(d,T_S(x))^2.
\end{equation}
\end{corollary}

\beginproof
Since the support function of every convex set is equal to the support function of its closure, without loss of generality, we may only consider nonempty closed convex subset $K(d)\subset T^{''}_K (g(x);\nabla g(x)d)$ or $\hat K(d)\subset T^2_K (g(x);\nabla g(x)d)$. It is easy to observe from Theorem \ref{swmc} that the optimal values of
\begin{equation}\label{minmax}
\inf \limits_{v\in \Theta}\sup \limits_{\lambda\in \Lambda^c(x;d)} -\lambda^Tv\quad
{\rm and}\quad \inf \limits_{w\in \Omega}\sup \limits_{\lambda\in \Lambda^c(x;d)} \nabla_{xx}^2L(x,\lambda)-\lambda^Tw-2\kappa(1-2\varepsilon)^2\|d\|^2
\end{equation}
are nonnegative, respectively, and that $\Lambda^c(x;d)$ is compact thanks to \cite[Lemma 6]{ye22}. Then by applying the minimax theorem \cite[Corollary 37.3.2]{rock} to programs \eqref{minmax}, we derive the assertions (i) and (ii) taking into account the fact that $-\sigma_{K(d)}$ and $-\sigma_{\hat K(d)}$ are upper semicontinuous.

For critical directions which are not necessarily $\varepsilon$-proximal normals, by applying the corresponding results of Theorem \ref{swmc}, we can also derive (i) and (ii) with inequality \eqref{nonprox2} in place of \eqref{prox2}.
\endproof

By \cite[Lemma 7]{ye22}, the directional nondegeneracy condition
\begin{equation}\label{dirnc}
\nabla g(\bar x)^T\lambda=0,\lambda\in {\rm span} (N_K(g(\bar x);\nabla g(\bar x)d))\Rightarrow \lambda=0
\end{equation}
at $\bar x$ in direction $d$ implies that the directional C- (M-) multiplier set $\Lambda(\bar x;d)=\Lambda^c(\bar x;d)$ is a singleton. It is easy to see that condition \eqref{dirnc} is stronger than the directional Robinson's CQ \eqref{dircq}.

The following result provides a sharper necessary optimality condition for second-order weak sharp minima of problem \eqref{op} under the directional nondegeneracy condition, which recovers \cite[Corollary 5]{ye22} if $\bar x$ reduces to be a local isolated optimal solution and $T_K^2(g(x);\nabla g(x)d) \neq \emptyset$. The proof is immediate from Theorem \ref{swmc} and hence it is omitted.

\begin{corollary}\label{swsm1}
Let $\bar x$ be a local second-order weak sharp minimizer of problem \eqref{op} with the corresponding constants $\kappa >0$ and $\delta >0$. Suppose that the directional nondegeneracy condition \eqref{dirnc} holds at every $x \in B_\delta (\bar x)$ in direction $d\in C(x)$. Then there is a unique multiplier $\lambda_0$ such that, for any $\varepsilon\in [0,1/2)$, $x\in S\cap B_{\delta}(\bar x)$, $d\in C(x)\cap N_S^{P,\varepsilon}(x)$, the following assertions hold:
\begin{enumerate}[{\rm (i)}]
\item $\sigma_{T_K^{''}(g(x);\nabla g(x)d)}(\lambda_0)= 0$;
\item If $T_K^2(g(x);\nabla g(x)d) \neq \emptyset$, then
\begin{equation}\label{prox3}
\nabla_{xx}^2L(x,\lambda_0)(d,d)-\sigma_{T_K^2(g(x);\nabla g(x)d)}(\lambda_0)\geq 2\kappa(1-2\varepsilon)^2\|d\|^2.
\end{equation}
\end{enumerate}

Besides, for any $x\in S\cap B_{\delta}(\bar x)$ and $d\in C(x)$, we also have (i) and (ii) hold true with \eqref{prox3} being replaced by
\begin{equation*}
\nabla_{xx}^2L(x,\lambda_0)(d,d)-\sigma_{T_K^2(g(x);\nabla g(x)d)}(\lambda_0)\geq 2\kappa\dist(d,T_S(x))^2.
\end{equation*}
\end{corollary}

%%%%%%%%%%%%%%%%%%%%%%%%%%%%%%%%%%%%%%%%%%%%%%%%%%%%%%%%%%%%%%%%%%%%%%%%%%%%%%%%%%%%%%%%%%%%%%%%%%%%%%%%%%%%%%%%%%%%%%%%%%%%%%%%%%%%%%%%%%%%%%%%%%%%%%%%%%%%%%%%%%%%%%%%%%%%%%%%%%%%%%%%%%%
\section{Sufficient conditions for second-order weak sharp minima}\label{3.2}

In this section, we aim to develop new sufficient second-order optimality conditions for second-order weak sharp minima of problem \eqref{op}.

To this end, we utilize a variation of the tangent cone which was first introduced in \cite{sward99} for studing weak sharp minima along with a generalized directional derivative that involves the solution set $S$. For a set $\Lambda \subset \mathbb{R}^n$ and $\bar x \in \Lambda$, the tangent cone to $\Lambda$ at $\bar x$ with respect to $S$ is defined by
$$T_{\Lambda,S}(\bar x):=\{d\in \mathbb{R}^n \mid \exists t_k\downarrow 0, u_k\xrightarrow[]{S}\bar x,d_k\rightarrow d \,\,{\rm with}\,\, u_k+t_kd_k\in \Lambda\}.$$

In general, we have $T_\Lambda(\bar x) \subset T_{\Lambda,S}(\bar x)$ and $T_\Lambda(\bar x) = T_{\Lambda,S}(\bar x)$ when $S=\{\bar x\}$.

The following result extends the fact that the tangent cone is a subset of a linearization cone to the tangent set with respect to the solution set.

\begin{proposition}\label{qc}
Given $\bar x\in S\subset \Phi:=g^{-1}(K)$ and $d\in T_{\Phi,S}(\bar x)$. Then
\begin{equation}\label{2cb0}
T_{\Phi,S}(\bar x)\subset \{d\in \mathbb{R}^n \mid \nabla g(\bar x)d\in T_{K,g(S)}(g(\bar x){)\}.}
\end{equation}
If, in addition, assume that MSCQ holds at $\bar x$ for the constraint system $g(x)\in K$ with modulus $\kappa$
and $g$ satisfying
\begin{equation}\label{dd}
g(x_k)\rightarrow g(\bar x)\Rightarrow x_k\rightarrow\bar x,
\end{equation}
then inclusion \eqref{2cb0} holds as equality and the following {estimation}
\begin{equation}\label{2cbb}
\dist(d,T_{\Phi,S}(\bar x))\leq \kappa \dist(\nabla g(\bar x)d, T_{K,g(S)}(g(\bar x)))
\end{equation}
{holds} for all $d\in \mathbb{R}^n$.
\end{proposition}
\beginproof
Let $d\in T_{\Phi,S}(\bar x)$, then there exist $t_k\downarrow 0, u_k\xrightarrow[]{S}\bar x,d_k\rightarrow d$ such that $u_k+t_kd_k\in \Phi$.
Then, we have
$$K\ni g(u_k+t_kd_k)=g(u_k)+t_k\nabla g(u_k)d_k+o(t_k).$$
Since $g$ is smooth, we have $g(u_k)\xrightarrow[]{g(S)}g(\bar x)$ and $\nabla g(u_k)d_k+\frac{o(t_k)}{t_k}\rightarrow \nabla g(\bar x)d$.
This shows that $\nabla g(\bar x)d \in T_{K,g(S)}(g(\bar x))$, and then \eqref{2cb0} holds true.

Next, assume that MSCQ holds at $\bar x$ for the constraint system $g(x)\in K$ with modulus $\kappa$.
Pick any $\varepsilon>0$ and $\kappa'>\kappa$. Let $d\in \mathbb{R}^n$ be fixed, then there exists $v \in T_{K,g(S)}(g(\bar x))$
such that
$\|\nabla g(\bar x)d-v\|<d(\nabla g(\bar x)d, T_{K,g(S)}(g(\bar x)))+\varepsilon.$
Since $v \in T_{K,g(S)}(g(\bar x))$, there exist $\{t_k\}\subset (0,+\infty)$, $\{u_k\} \subset S$ and $\{v_k\}\subset \mathbb{R}^m$ such that
$t_k\rightarrow 0, g(u_k)\rightarrow g(\bar x), v_k\rightarrow v$ and $g(u_k) +t_k v_k\in K$.
By \eqref{dd}, we have $u_k\rightarrow \bar x$.
Let $x_k:=u_k +t_k d$, then, for sufficiently large $k$, it follows from the  MSCQ that
\begin{equation*}
\aligned
d(u_k +t_k d,\Phi)&\leq \kappa' d(g(u_k +t_k d),K)\\
&\leq \kappa' \|g(u_k +t_k d)-g(u_k)-t_k v_k\|\\
&\leq \kappa' t_k\left\|\nabla g(u_k)d-v_k\right\|+o(t_k).
\endaligned
\end{equation*}
Then, there exists $x_k'\in \Phi$ such that
\begin{equation}\label{xkp1}
\left\|\frac{x_k'-u_k}{t_k}-d\right\|\leq \left\|\nabla g(u_k)d-v_k\right\|+\frac{o(t_k)}{t_k}.
\end{equation}
Since the right side of \eqref{xkp1} converges to $\kappa'\left\|\nabla g(\bar x)d-v\right\|$, we know that $\{\frac{x_k'-u_k}{t_k}\}$ is bounded. Without loss of generality, we assume that $\frac{x_k'-u_k}{t_k}\rightarrow d'$, then $d'\in T_{\Phi,S}^{''}(\bar x)$.
Passing to the limit in \eqref{xkp1}, we arrive at
$$\dist(d,T_{\Phi,S}(\bar x))\leq\|d-d'\|\leq \kappa'\left\|\nabla g(\bar x)d-v\right\|\leq\kappa'\dist(\nabla g(\bar x)d, T_{K,g(S)}(g(\bar x)))+\varepsilon.$$
Because $\kappa'$ can be chosen arbitrarily close to $\kappa$ and $\varepsilon$ can be chosen  arbitrarily small, we know that \eqref{2cbb} holds true.
 From \eqref{2cbb}, we may conclude that
$$\{d\in \mathbb{R}^n \mid \nabla g(\bar x)d\in T_{K,g(S)}(g(\bar x))\}\subset T_{\Phi,S}(\bar x), $$
which indicates that inclusion \eqref{2cb0} holds as equality.
\endproof

Exploiting the spirit of the above definition, we introduce the following variations of the second-order tangent sets which take into account of the level set $S$.

\begin{definition}[Second-order Tangent Sets w.r.t. Level Sets]\label{t1s}
Given sets $\Lambda,S\subset \mathbb{R}^n$, $\bar x\in  \Lambda$ and $d \in T_{\Lambda,S}(\bar x)$.

(i) The outer second-order tangent set to $\Lambda$ at $\bar x$ in direction $d$ with respect to $S$ is defined by
\begin{equation*}
\aligned
T_{\Lambda,S}^2(\bar x;d):=&\{w\in \mathbb{R}^n \mid \exists t_k\downarrow 0,u_k\xrightarrow[]{S}\bar x,w_k\rightarrow w\,\,{\rm such\,that}\\
&u_k +t_k d+\frac{1}{2}t_k^2w_k\in \Lambda\}.
\endaligned
\end{equation*}

(ii) The asymptotic second-order tangent cone to $\Lambda$ at $\bar x$ in direction $d$ with respect to $S$ is defined by
\begin{equation*}
\aligned
T_{\Lambda,S}^{''}(\bar x;d):=&\{w\in \mathbb{R}^n \mid \exists (t_k,r_k)\downarrow (0,0),u_k\xrightarrow[]{S}\bar x, w_k\rightarrow w\,\,{\rm such\,that}\\
&t_k/r_k\rightarrow 0,u_k +t_k d+\frac{1}{2}t_kr_kw_k\in \Lambda\}.
\endaligned
\end{equation*}
\end{definition}

It is clear that if $d\not\in T_{\Lambda,S}(\bar x)$, we have $T_{\Lambda,S}^{2}(\bar x;d)=T_{\Lambda,S}^{''}(\bar x;d)=\emptyset$. Observe also that if $S = \{\bar x\}$, we have that $T_{\Lambda}^{2}(\bar x;d)=T_{\Lambda,S}^{2}(\bar x;d)$ and $T_{\Lambda}^{''}(\bar x;d)=T_{\Lambda,S}^{''}(\bar x;d)$ in this case. Furthermore, we can see from the proof of Theorem \ref{sff} that, for any $d \in T_{\Lambda,S}(\bar x)$,  both $T_{\Lambda,S}^{2}(\bar x;d)$ and $T_{\Lambda,S}^{''}(\bar x;d){\setminus \{0\}}$ cannot be empty simultaneously.

\begin{proposition}\label{qcb}
Given $\bar x\in S\subset \Phi:=g^{-1}(K)$ and $d\in T_{\Phi,S}(\bar x)$. Then, for any $\varepsilon>0$,
\begin{equation}\label{bc}
T_{\Phi,S}^2(\bar x;d)\subset \{w\in \mathbb{R}^n \mid \nabla g(\bar x)w+\nabla^2 g(\bar x)(d,d)\in T_{K,g_\varepsilon(S)}^2(g(\bar x);\nabla g(\bar x)d)\}
\end{equation}
and
\begin{equation}\label{bc1}
T_{\Phi,S}^{''}(\bar x;d)\subset \{w\in \mathbb{R}^n \mid \nabla g(\bar x)w\in T_{K,g_\varepsilon(S)}^{''}(g(\bar x);\nabla g(\bar x)d)\},
\end{equation}
where $g_\varepsilon(S)=g(S)+\varepsilon {\rm diam}(S)B_{\mathbb{R}^m}$.
\end{proposition}

\beginproof
Pick any $\varepsilon>0, w\in T_{\Phi,S}^2(\bar x;d)$, then there exist $t_k\downarrow 0, u_k\xrightarrow[]{S}\bar x$ and $w_k\rightarrow w$ such that
$x_k:=u_k+t_k d+\frac{1}{2}t_kw_k\in \Phi$ for all $k\in \mathbb{N}$.
Then,
\begin{equation}\label{bhh}
\aligned
K&\ni g(x_k)=g(u_k)+t_k\nabla g(u_k)d+\frac{1}{2}t_k^2\nabla g(u_k)w_k+\frac{1}{2}t_k^2(\nabla g(u_k)w_k+\nabla^2 g(\bar x)(d,d))+o(t_k^2)\\
&=g(u_k)+t_k(\nabla g(u_k)-\nabla g(\bar x))d+t_k\nabla g(\bar x)d+\frac{1}{2}t_k^2(\nabla g(u_k)w_k+\nabla^2 g(\bar x)(d,d))+o(t_k^2).
\endaligned
\end{equation}
Since $g$ is twice continuously differentiable, then there exists $l>0$, such that $\|\nabla g(u_k)-\nabla g(\bar x)\|\leq l\|u_k-\bar x\|$.
By $t_k\downarrow 0$ and $u_k\xrightarrow[]{S}\bar x$, we have $\|t_k(\nabla g(u_k)-\nabla g(\bar x))d\|<\varepsilon {\rm diam}(S)$ for sufficiently large $k$. And then, we have
$g(u_k)+t_k(\nabla g(u_k)-\nabla g(\bar x))d\xrightarrow[]{g_\varepsilon(S)}g(\bar x)$.
Note that
$$\nabla g(u_k)w_k+\nabla^2 g(\bar x)(d,d)+\frac{o(t_k^2)}{\frac{1}{2}t_k^2}\to \nabla g(\bar x)w+\nabla^2 g(\bar x)(d,d),$$
then, \eqref{bhh} implies $\nabla g(\bar x)w+\nabla^2 g(\bar x)(d,d)\in T_{K,g_\varepsilon(S)}^2(g(\bar x);\nabla g(\bar x)d)$, and then \eqref{bc} holds true.

For \eqref{bc1}, we can { prove it} similarly, so we omit it.
\endproof
By Proposition \ref{qcb}, one has,  for any $\varepsilon>0$,
\begin{eqnarray}
\Omega_S &:=& \nabla g(\bar x)(T_{\Phi, S}^2(\bar x;d))+\nabla^2 g(\bar x)(d,d)\subset T_{K,g_\varepsilon(S)}^2(g(\bar x);\nabla g(\bar x)d),\label{strin}\\
\Theta_S  &:=& \nabla g(\bar x)(T_{\Phi, S}^{''}(\bar x;d))\subset T_{K,g_\varepsilon(S)}^{''}(g(\bar x);\nabla g(\bar x)d).\label{strintheta}
\end{eqnarray}

Following the aforementioned definition, we introduce below a variation of the critical cone with respect to the level set $S$:
$$C_S(\bar x):=\{d\in T_{\Phi,S}(\bar x)\mid\nabla f(x)d\leq 0\ {\rm for\ all}\  x\in {\rm bd}(S)\ {\rm close\ to} \ \bar x\}.$$
It is clear that $C_S(\bar x)\subset T_{\Phi,S}(\bar x)$ and $C_S(\bar x)=C(\bar x)$ for the case of $S=\{\bar x\}$.

{Next, using the previously introduced second-order tangent sets with respect to the solution set, we give a point-based second-order sufficient optimality condition for local weak sharp minima of problem \eqref{op}. }

\begin{theorem}[Point-based Second-order Sufficient Condition]\label{sff}
Assume that for every $d\in T_{\Phi,S}(\bar x) \cap N_S(\bar x)\backslash\{0\}$, we have $\nabla f(x)d=0$ for all $x\in {\rm bd}(S)$ close to $\bar x$ and there is some $\lambda\in \mathbb{R}^m$ such that
$\nabla_x L(\bar x,\lambda)=0$ and the following conditions hold with some $\kappa>0$:
\begin{enumerate}[{\rm (i)}]
\item $\langle\lambda,v\rangle< 0,\quad\forall v\in \nabla g(\bar x)(T_{\Phi,S}^{''}(\bar x,d)\cap\{d\}^\bot\backslash\{0\})$;
\item $\nabla_{xx}^2L(\bar x,\lambda)(d,d)-\sigma_{\nabla g(\bar x)(T_{\Phi,S}^2(\bar x;d)\cap\{d\}^\bot)+\nabla^2 g(\bar x)(d,d)}(\lambda)>\kappa\|d\|^2$.
\end{enumerate}
Then, there exists $\delta>0$ such that
\begin{equation}\label{2wk}
f(x)-f(\bar x)\geq \kappa[\dist(x,S)]^2 \quad \forall x\in B_\delta(\bar x).
\end{equation}
Moreover,
$\bar x$ is a local second-order weak sharp minimizer of problem \eqref{op}.
\end{theorem}

\beginproof
Assume to the opposite that \eqref{2wk} does not hold, then there exists a sequence
$\{x_k\}\subset \Phi\cap B_{1/k}(\bar x)$ such that
\begin{equation}\label{fc}
f(x_k)<f(\bar x)+\kappa[\dist(x_k,S)]^2.
\end{equation}
Pick $u_k\in S$ such that $\|x_k-u_k\|=\dist(x_k,S)$.
Let $t_k=\|x_k-u_k\|$ and $d_k=\frac{x_k-u_k}{t_k}\in S_{\mathbb{R}^n}$.
Since $t_k=\dist(x_k,S)\leq \|x_k-\bar x\|\to 0$, we have $t_k\downarrow 0$ and $u_k\rightarrow \bar x$.
Without loss of generality, we may assume that $d_k$ converges to some vector $d \in T_{\Phi,S}(\bar x)\cap S_{\mathbb{R}^n}$.
Since $u_k\in P_S(x_k)$, we have $d_k\in N_S^P(u_k)$, and then $d\in N_S(\bar x)$.
Now we consider the sequence $w_k:=\frac{x_k-u_k-t_kd}{\frac{1}{2}t_k^2}$ and arrive at the following two situations:

(a) If $\{w_k\}$ is bounded, then we have $w_k=\frac{2(d_k-d)}{t_k}\rightarrow w$
for some $w\in T_{\Phi,S}^2(\bar x;d)$ (take a subsequence if necessary).

(b) If $\{w_k\}$ is not bounded, then we may assume that $s_k=\|w_k\|\rightarrow \infty$ and $\tilde{w}_k:=\frac{w_k}{\|w_k\|}\rightarrow w\in S_{\mathbb{R}^n}$ (take a subsequence if necessary). Note that
\begin{equation*}
\tilde{w}_k:=\frac{w_k}{\|w_k\|}=\frac{x_k-\tilde{x}_k-t_kd}{\frac{1}{2}t_kr_k}=\frac{2(d_k-d)}{r_k}\rightarrow w,
\end{equation*}
where $r_k:=t_ks_k$, $r_k=\|\frac{x_k-u_k-t_kd}{\frac{1}{2}t_k}\|\rightarrow 0$ and $t_k/r_k=1/s_k\rightarrow 0$. Then we have $w\in T_{\Phi,S}^{''}(\bar x;d)\backslash\{0\}$.

Furthermore, we claim that $w\in \{d\}^\bot$ for both cases (a) and (b). Indeed, since $d_k\in S_{\mathbb{R}^n}$, we deduce from the above expressions for $w_k$ and $\tilde{w}_k$ that
$w\in T_{S_{\mathbb{R}^n}}(d) = \{d\}^\bot$. By assumptions, we take $\lambda\in \mathbb{R}^m$ such that $\nabla_x L(\bar x,\lambda)=0$ and conditions (i) and (ii) hold. Then, $\nabla f(\bar x)+\nabla g(\bar x)^T\lambda=0$.

For case (a), we have $x_k=u_k+t_kd+\frac{1}{2}t_k^2w_k$.
Note that $\tilde{x}_k\xrightarrow[]{{\rm bd}(S)}\bar x$, one has $\nabla f(u_k)d=0$ for sufficiently large $k$ and $f(u_k)=f(\bar x)$. It follows from \eqref{fc} that
$$f(x_k)-f(u_k)= t_k\nabla f(u_k)d+\frac{1}{2}t_k^2\left(\nabla f(u_k)w_k+\nabla^2 f(u_k)\left(d,d\right)\right)+o(t_k^2)<\kappa t_k^2.$$
Dividing both sides of above inequality by $\frac{1}{2}t_k^2$ and passing to the limit, we have
\begin{equation}\label{ww}
\nabla f(\bar x)w+\nabla^2 f(\bar x)\left(d,d\right)\leq \kappa.
\end{equation}
For convenience, let $\Lambda:=\nabla g(\bar x)(T_{\Phi,S}^2(\bar x;d)\cap\{d\}^\bot)+\nabla^2 g(\bar x)(d,d)$. Since $w\in T_{\Phi,S}^2(\bar x;d)\cap\{d\}^\bot$, we have
$\nabla g(\bar x)w+\nabla^2 g(\bar x)\left(d,d\right)\in \Lambda.$
And then
\begin{equation}\label{ww1}
\left\langle\lambda,\nabla g(\bar x)w+\nabla^2 g(\bar x)\left(d,d\right)\right\rangle
\leq \sigma_{\Lambda}(\lambda).
\end{equation}
Adding \eqref{ww} and \eqref{ww1}, we obtain that
\begin{equation*}
\aligned
\sigma_{\Lambda}(\lambda)+\kappa&\geq\nabla f(\bar x)w+\langle\nabla g(\bar x)^T\lambda,w\rangle+\nabla^2 f(\bar x)\left(d,d\right)+
\left\langle\lambda,\nabla^2 g(\bar x)\left(d,d\right)\right\rangle\\
&=0+\nabla^2_{xx} L(\bar x,\lambda)\left(d,d\right)= \nabla^2_{xx} L(\bar x,\lambda)\left(d,d\right),
\endaligned
\end{equation*}
which is a contradiction to assumption (ii).

For case (b), we have $x_k=u_k+t_kd+\frac{1}{2}t_kr_k\tilde{w}_k$.
It follows from \eqref{fc} that
$$f(x_k)-f(u_k)= t_k\nabla f(u_k)d+\frac{1}{2}t_kr_k\nabla f(u_k)\tilde{w}_k+\frac{1}{2}t_k^2\nabla^2 f(u_k)\left(d,d\right)+o(t_k^2)<\frac{1}{k}t_k^2.$$
Since $t_k/r_k\rightarrow 0$, dividing both sides of above inequality by $\frac{1}{2}t_kr_k$ and passing to the limit, we obtain that
$-\langle \lambda, \nabla g(\bar x)w\rangle=\nabla f(\bar x)w\leq 0,$ which is a contradiction to (i). Therefore, we conclude that $\bar x$ is a local second-order weak sharp minimizer of problem \eqref{op}.
\endproof
{ It is clear that when $\bar x$ is an isolated point of $S$, Theorem \ref{sff} reduces to \cite[Theorem 4.14]{OuyangYe}.}
It follows from \eqref{strin}, \eqref{strintheta} and Theorem \ref{sff} that  we have following explicit sufficient condition.

 \begin{corollary}\label{Cor4.5}
%Let $S$ be a closed subset of feasible set $\Phi$  such that $f$ is constant on it, $\bar x \in S$ and $\kappa>0$. 
Assume that for every nonzero vector $d\in T_{\Phi,S}(\bar x) \cap N_S(\bar x)$, we have $\nabla f(x)d=0$ for all $x\in {\rm bd}(S)$ close to $\bar x$ and there is some $\lambda\in \mathbb{R}^m$ and $\varepsilon>0$ such that
$\nabla_x L(\bar x,\lambda)=0$ and the following conditions hold:
\begin{enumerate}[{\rm (i)}]
\item $\langle\lambda,v\rangle< 0,\quad\forall v\in T_{K,g_\varepsilon(S)}^{''}(g(\bar x); \nabla g(\bar x)d)\cap \nabla g(\bar x)(\{d\}^\bot\backslash\{0\})$;
\item $\nabla_{xx}^2L(\bar x,\lambda)(d,d)-\sigma_{T_{K,g_\varepsilon(S)}^2(g(\bar x); g(\bar x)d)}(\lambda)>\kappa\|d\|^2$.
\end{enumerate}
Then, there exists $\delta>0$ such that
\begin{equation*}
f(x)-f(\bar x)\geq \kappa[\dist(x,S)]^2 \quad \forall x\in B_\delta(\bar x).
\end{equation*}
Moreover,
$\bar x$ is a local second-order weak sharp minimizer of problem \eqref{op}.
\end{corollary}

Since a local minimizer satisfying the second-order growth condition is always a local second-order weak sharp minimizer (as outlined in conditions \eqref{los} and \eqref{quadraticgrow}), we derive the following sufficient condition.
\begin{theorem}\label{sff3} 
 Let $\bar x $ be an isolated point of $S$. Suppose that $\nabla f(\bar x)d \geq 0$ for all $d \in T_{\Phi}(\bar x)$ and there exists {$\lambda\in \mathbb{R}^m$ with $\nabla_x L( \bar x, \lambda)=0$} such that for any $d \in C(\bar x) \setminus \{0\}$,
\begin{enumerate}[{\rm (i)}]
\item $\langle \lambda,v \rangle <0, \quad \forall v \in \nabla g(\bar x)(T_\Phi^{''}(\bar x,d)\cap (\{d\}^\perp \setminus \{0\}))$;
\item $\nabla_{xx}^2L(\bar x,\lambda)(d,d)-\sigma_{\nabla g(\bar x)(T_\Phi^2(\bar x;d)\cap \{d\}^\perp)+\nabla^2 g(\bar x)(d,d)}(\lambda)> 0$.
\end{enumerate}
Then  $\bar x$ is a local second-order weak sharp minimizer of problem \eqref{op}.
\end{theorem}
\beginproof
Since $\bar x \in S$, we have $\|x-\bar x\| \geq \dist(x,S)$ for any $x$ around $\bar x$. Thus, the result follows from {By Theorem 4.14 of \cite{OuyangYe}.
\endproof

}

%The sufficient conditions in Theorem \ref{sff} are formulated using second-order tangent sets/asymptotic second-order tangent cones relative to level sets, which are generally larger than standard second-order tangent sets/asymptotic second-order tangent cones. Consequently, compared to Theorem \ref{sff}, the sufficient conditions in Theorem \ref{sff3} are weaker. However, they lead to a special type of local second-order weak sharp minimizer, which is an isolated point in the solution set $S$.

The sufficient conditions in Theorem \ref{sff} are formulated using second-order tangent sets/asymptotic second-order tangent cones relative to level sets, which are generally larger than standard second-order tangent sets/asymptotic second-order tangent cones used in Theorem \ref{sff3}.  {On the other hand, compared to Theorem \ref{sff3}, Theorem \ref{sff} refines the range of directions $d$ by utilizing the normal cone of the solution set $S$, making it easier to verify. In fact, Theorem \ref{sff3} leads to a special type of local second-order weak sharp minimizer, which is an isolated point in the solution set $S$.}

{
\section{Conclusions}
In this paper, we establish both neighborhood-based necessary conditions and point-based sufficient conditions for second-order weak sharp minima in general nonconvex set-constrained optimization problems. These results are obtained without assuming the convexity of the constraint set or the second-order tangent set. By using generalized support functions and directional constraint qualifications, our approach avoids reliance on uniform regularity assumptions. While our results broaden the scope of second-order theory for weak sharp minima, challenges remain in computing the second-order tangent sets w.r.t. level sets, which arise in the sufficient optimality conditions. We leave this for future research.
}

%\section*{Acknowledgments}
%This work was partially supported by the China
%Scholarship Fund (202208535016). Wei Ouyang would like to thank Professor Jane J. Ye for her warm hospitality during her visit at the Department of Mathematics and Statistics,  University of Victoria.

%%%%%%%%%%%%%%%%%%%%%%%%%%%%%%%%%%%%%%%

\end{sloppypar}
\end{document}